\documentclass[article]{amsart}

\usepackage{CJK,CJKnumb}
\usepackage{color}
\usepackage{indentfirst}
\usepackage{latexsym,bm}
\usepackage{amsmath,amssymb}
\usepackage{graphicx}
\usepackage{bbm}
\usepackage{amsthm}
\usepackage[ansinew]{inputenc}
\usepackage{array}
\usepackage{amsxtra}
\usepackage{amstext}
\usepackage{latexsym}
\usepackage{amsfonts}
\usepackage{bm}
\usepackage[pdftex,pdfstartview=FitH,bookmarksnumbered,bookmarksopen,bookmarksopenlevel=2,CJKbookmarks]{hyperref}

\newtheorem{theorem}{Theorem}[section]
\newtheorem{definition}[theorem]{Definition}

\newcommand{\secref}[1]{Section \ref{#1}}

 \numberwithin{equation}{section}

\begin{document}

{\title
[
From conservative to dissipative systems
]
{{
From conservative to dissipative systems
through quadratic change of time, with application to the curve-shortening flow}}
}
\author{
Yann Brenier and Xianglong Duan}
\date\today



\address{
CNRS UMR 7640 \\ 
Ecole Polytechnique   \\ 
Palaiseau\\
France}

\email{yann.brenier@polytechnique.edu and xianglong.duan@polytechnique.edu}

\begin{abstract}

We provide several examples of dissipative systems that can be obtained from conservative ones through a simple, quadratic,
change of time. A typical example is the curve-shortening flow in $\mathbb{R}^d$, which is a particular case of
mean-curvature flow with co-dimension higher than one (except in the case $d=2$).
Through such a change of time, this flow can be formally derived from the conservative model of vibrating strings obtained
from the Nambu-Goto action. Using the concept of ``relative entropy" (or ``modulated energy"), borrowed from the theory
of hyperbolic systems of conservation laws, we introduce a notion of generalized solutions,
that we call dissipative solutions, for the curve-shortening flow.  For given initial conditions, the set of generalized solutions
is convex, compact, if not empty. Smooth solutions to the curve-shortening flow are always unique in this setting.

\end{abstract}
\maketitle

\section*{Introduction}

There are many examples of dissipative systems that can be derived from conservative ones.
A classical example is the heat equation (or more generally the so-called ``porous medium"
equation) that can be derived from the Euler equations of isentropic gases. The derivation can be
done in many different ways, for example by adding a very strong friction term or by homogenization techniques
or by properly rescaling the time variable by a small
parameter (through the so-called ``parabolic scaling"). In the present paper, we will focus on a very straightforward
idea (that does not seem to be popular, to the best of our knowledge): just perform the
quadratic change of time
$t\rightarrow \theta=t^2/2$.
In \secref{examples},
we provide several examples: we start with the very simple example of conservative forces in classical mechanics (with the
Galileo model of falling bodies as a borderline case).
Next, we briefly retrieve from the Euler equation of isentropic gases
the Darcy law and the porous
medium equation, and, in particular the heat equation from the Euler equation of isothermal gases.
Our third example, at the interface of Geometry and High Energy Physics, starts with the
conservative evolution of classical strings according to the Nambu-Goto action,
from which we get, by quadratic change of time,
the dissipative geometric model of curve-shortening in $\mathbb{R}^d$,
which is the simplest example of mean-curvature flow with co-dimension higher than $1$:
\begin{equation}
\label{shortening}
\partial_\theta X
=\frac{1}{|\partial_s X|}\partial_s(\frac{\partial_s X}{|\partial_s X|}),
\end{equation}
where $s\rightarrow X(\theta,s)$ describes a time-dependent curve in $\mathbb{R}^d$ and $|\cdot|$ denotes
the Euclidean norm.
\\
\\
In \secref{main}, we will finally discuss the system of PDEs
$$
\partial_\theta B+\nabla\cdot \left(\frac{B\otimes P-P\otimes B}{\rho}\right)=0,\;\;\;\nabla\cdot B=0,
$$
$$
P=\nabla\cdot\left(\frac{B\otimes B}{\rho}\right),
\;\;\;\rho=|B|
$$
which turns out to be nothing but the ``Eulerian version" (in $\mathbb{R}^d$) of the curve-shortening
model (\ref{shortening}). Typically, in the case of a single loop $X$ subject to the curve-shortening flow,
$B$ would just be the singular
vector-valued measure
$$
(\theta,x)\rightarrow B(\theta,x)=\int_{\mathbb{R/Z}}\delta(x-X(\theta,s))\partial_sX(\theta,s) ds\in\mathbb{R}^d,
$$
for which the system of PDE makes sense since all nonlinearities are homogeneous of degree one.
(See Appendix $1$ for more details.)
These equations admit a ``non-conservative" version
\begin{equation}\label{eq:shortnc}
\partial_\theta b+(v\cdot\nabla)b=(b\cdot\nabla)v+b v^2,
\;\;\;v=(b\cdot\nabla)b,
\end{equation}
for the reduced variables $b=B/|B|$ and $v=P/|B|$.
For the conservative system,
we define a concept of ``dissipative solutions"
related to the work of P.-L. Lions for the Euler equation of incompressible fluids \cite{Li} (see also \cite{BDS,Vo}) or to
the work of L. Ambrosio, N. Gigli, G. Savar\'e \cite{AGS} for the heat equation
and, overall, quite similar to the one recently introduced by the first author in \cite{Br-topo}. We also refer to the works of
A. Tzavaras and collaborators \cite{DeStTz,LT}, E. Feireisl and collaborators \cite{Fe,FGSW} for related concepts of ``dissipative solutions".
The main point of the present paper is to show how to get the formulation right:
we start from the Eulerian version of the string equation, for which we can
use the ``relative entropy" method, quite classical in the theory of hyperbolic systems of conservation laws to get ``weak-strong" uniqueness
results (see \cite{Da}),
then, we apply the quadratic change of
time to get a good concept
of dissipative solutions for the curve-shortening flow, namely:

\begin{definition}\label{def:dddf1}
Let us fix $T>0$ and denote $\mathbb{T}^d=\mathbb{(R/Z)}^d$.
We say that $(B,P)$ with
$$
B\in C([0,T],C(\mathbb{T}^d
,\mathbb{R}^{d})'_{w^{*}}),\;\;\;P\in C([0,T]\times \mathbb{T}^d,\mathbb{R}^{d})'
$$
is a dissipative solution of the curve-shortening flow
with initial data $B_{0}\in C(\mathbb{T}^d,\mathbb{R}^{d})'$ if and only if:
\\
\\
i)
$\;\;\;B(0)=B_0$,
$\;\;\nabla\cdot B=0$
in sense of distributions;
\\
\\
ii) $B$ and $P$ are bounded, respectively in the spaces
$C^{1/2}([0,T],(C^{1}(\mathbb{T}^d))'_{w^*})$ and
$C([0,T]\times \mathbb{T}^d,\mathbb{R}^{d})',$
by constants depending only on $T$ and $\int_{\mathbb{T}^d}|B_0|$.
\\
\\
iii) For all $\lambda>0$, $\theta\in [0,T]$,
for all smooth trial functions $(b^*,v^*,A)$ valued in $\mathbb{R}^d$,
with $\|A\|_{\infty}\le \lambda$ and ${b^*}^2=1$, for all $r\ge c^*+\frac{\lambda^2}{2}+\lambda\|v^*\|_{\infty}$,
where $c^*$ is a constant depending
explicitly on $(b^*,v^*)$, we have:
\begin{multline}\label{dissipative-1}
e^{-r\theta}
\int \eta(\theta)
+\int_0^\theta e^{-r\sigma}
\left[\int P\cdot (A-L_3)+\left(r-c^*-\frac{A\cdot(A+2v^*)}{2}\right)\eta\right.\\
\left.-B\cdot \left(L_2+b^*\frac{A\cdot(A+2v^*)}{2}\right)\right]
(\sigma)d\sigma\le \int \eta(0).
\end{multline}
where
\begin{equation}\label{dissipative-4}
\eta=|B|-B\cdot b^*,
\end{equation}
\begin{equation}\label{dissipative-2}
\mathrm{L}_2=-\partial_\theta b^*-(v^*\cdot\nabla)b^*+(b^*\cdot\nabla)v^*+b^*{v^*}^2
-b^*(b^*\cdot\nabla)(b^*\cdot v^*),
\end{equation}
\begin{equation}\label{dissipative-3}
\mathrm{L}_3=-v^*+(b^*\cdot\nabla)b^*.
\end{equation}
\end{definition}
Here $C(\mathbb{T}^3,\mathbb{R}^d)'_{w^{*}}$ is metrizable space, we can equip a metric that is consistent with the weak-* topology. The ``weak compactness" of such solutions (i.e. any sequence of dissipative solutions has accumulations
points, in a suitable weak sense, and each of them is still a dissipative solution) directly follows from:

\begin{theorem}\label{theorem1}
For fixed initial condition $B_0$, the set of dissipative solutions, if not empty, is convex and compact
for the weak-* topology of measures.
\end{theorem}

Notice that it is more challenging to prove that the set of dissipative solutions is not empty. The standard strategy is as follows:
\\
i) construct
smooth approximate solutions $(B^\epsilon,P^\epsilon)$ with smooth approximate initial data $B^\epsilon_0$;
\\
ii) show that, the approximate solutions are relatively compact for the weak-* topology of measures, and, for any trial functions $(b^*,v^*,A)$, satisfy inequalities \eqref{dissipative-1} with some small error terms;
\\
iii) let $\epsilon$ go to zero, and prove that the limit $(B,P)$ is a dissipative solution.
\\
\\
To keep this paper simple, we leave this (important) step for a future work
\cite{XD}, in the more
general framework of the Born-Infeld theory \cite{BI,Br-BI}.
Finally,
we establish a ``weak-strong" uniqueness principle in the following sense:

\begin{theorem}\label{theorem2}

Let $(b,v)\in C^{1}([0,T]\times\mathbb{T}^d;\mathbb{R}^d\times \mathbb{R}^d)$,
be a smooth solution of the non-conservative
form of the curve-shortening flow (\ref{eq:shortnc}) with $b^2=1$.
Then any dissipative solution satisfies $B=|B| b$ and $P=|B| v,$ as soon as
$B(0)=|B(0)|b(0)$.
\end{theorem}

Notice that this is not a full uniqueness result: only the homogeneous variables $b=B/|B|$ and
$v=P/|B|$ get unique and a lot of room is left for the evolution of $|B|$ itself. So the concept of
dissipative solutions seems to suffer from the same type of ambiguity as the more general concept of Brakke solutions for mean-curvature flows
\cite{Bra}.
\subsection*{Acknowledgements}
This work has been partly supported by the contract : MAGA ANR-16-CE40-0014.  
\\
The first author would like to thank
the Erwin Schr\"odinger Institute (ESI) for its hospitality when this work was started.
He is also very grateful to Dmitry Vorotnikov for explaining to him, at ESI in the summer of 2016, the possibility of deriving some mean-curvature motions
as gradient flows in optimal transportation style, in the spirit of \cite{JKO,Ot,Sa,Vi}. This information was very useful to understand that
the methods used for the Born-Infeld equations by the second author \cite{XD} could also be applied to the curve-shortening flow.

\section{Examples of quadratic change of time}
\label{examples}
\subsection
{Quadratic change of time of a simple dynamical system
}

Under the quadratic change of time $t\rightarrow \theta=t^2/2$, the ordinary dynamical system
$$
\frac{d^2 X}{dt^2}=-\nabla \varphi(X)
$$
becomes
$$
-\nabla \varphi(X)=\frac{d}{dt}(\frac{dX}{d\theta}\frac{d\theta}{dt})
=\frac{d}{dt}(t\frac{dX}{d\theta})
=\frac{dX}{d\theta}+
t\frac{d\theta}{dt}\frac{d^2 X}{d\theta^2}
=\frac{dX}{d\theta}+
2\theta\frac{d^2 X}{d\theta^2}
$$
with two asymptotic regimes as $\theta$ becomes either very small or very large:
\\
the ``gradient flow"
$$
\frac{dX}{d\theta}
=-\nabla \varphi(X),
$$
and the inertial motion
$$
\frac{d^2X}{d\theta^2}=0.
$$
Notice, in the first case, that only the initial position can be chosen freely, since
$$
\frac{dX}{dt}=\frac{dX}{d\theta}\frac{d\theta}{dt}
=t\frac{dX}{d\theta}
$$
necessarily vanishes at $t=0$.
Consistently, the  conservation of energy in the original time variable reads
$$
\frac{d}{d t}(\frac{1}{2}|\frac{dX}{d t}|^2+\varphi(X))=0
$$
and becomes, with the new time variable $\theta=t^2/2$
$$
\frac{d}{d\theta}[\varphi(X)]+
\theta
\frac{d}{d\theta}|\frac{dX}{d\theta}|^2
=-|\frac{dX}{d\theta}|^2,
$$
leading
to the dissipation of energy
$$
\frac{d}{d\theta}[\varphi(X)]
=-|\frac{dX}{d\theta}|^2,
$$
 in the asymptotic gradient flow regime.
Furthermore, we may compare the respective solutions $X(t)$ and $Z(\theta)$
of the dynamical system and the gradient flow,
with initial conditions
$$
X(t=0)=X_0=Z(\theta=0)
,\;\;\;\frac{dX}{dt}(t=0)=0,
$$
just by
monitoring
the ``modulated energy" (or ``relative entropy")
\begin{equation}\label{relative}
\frac{1}{2}|\frac{dX}{dt}-t\frac{dZ}{d\theta}|^2
+\varphi(X)-\varphi(Z)-\nabla \varphi(Z)\cdot(X-Z),
\end{equation}
provided $\varphi$ is strongly convex with bounded third derivatives. We get, after
elementary calculations,
\begin{equation}
\label{error}
|X(t)-Z(t^2/2)|^2+|\frac{dX}{dt}(t)-t\frac{dZ}{d\theta}(t^2/2)|^2\le Ct^5 ,\;\;\;\forall t\in [0,T],
\end{equation}
where $C$ is a constant that depends only on $T$, $Z$ and potential $\varphi$.
(Notice that the smallest expected error is $O(t^6)$ as shown by the example
$d=1$, $\varphi(x)=|x|^2/2$, for which $X(t)=X(0)\cos(t)$, while $Z(\theta)=X(0)\exp(-\theta)$.)
More details on the concept of "modulated energy" and the proof of (\ref{error}) can be found in an appendix at the
end of this paper.

\subsection*{Remark: the Galileo experiment}

The quadratic change of time $t\rightarrow \theta=t^2/2$ remarkably fits with the famous
experiment by Galileo, which was the starting point of modern classical mechanics:
a rigid ball descends a rigid ramp of constant slope,
with zero initial velocity and constant acceleration $G$, reaching position
$$
X=x_0+\frac{Gt^2}{2}
$$
at time $t$.
So, $X$ is just a linear function of the rescaled time $\theta$,
$X=x_0+\theta G$ and we not only get
$$
\frac{dX}{d\theta}+
{{2\theta\frac{d^2 X}{d\theta^2}}}=G
$$
but also $simultaneously$
$$
\frac{dX}{d\theta}=G,\;\;\;
{{\frac{d^2 X}{d\theta^2}}}=0,
$$
i.e. $both$ gradient flow $and$ inertial motion, with respect to the rescaled time $\theta$.

\subsection{From the Euler equations to the heat equation and the Darcy law}

Let us now move to a PDE example and explain how the
Darcy law and the ``porous medium" equation (and, in particular, the standard heat equation) can be recovered by quadratic change of time from
the Euler equations of isentropic compressible fluids. These equations read
\begin{equation}\label{}
\partial_t\rho+\nabla\cdot(\rho v)=0,\;\;\;
\partial_t (\rho v)+\nabla\cdot(\rho v\otimes v)=-\nabla  p,
\end{equation}
where ${(\rho,p,v)\in\mathbb {R}^{1+1+3}}$
are the density, pressure and velocity fields of the fluid, $p$ being a given function
of $\rho$ (such as $p=\rho$, in the ``isothermal" case).
We set
\begin{equation}\label{}
t\rightarrow \theta=t^2/2,
\;\;\;\;
\rho(t,x)\rightarrow \rho(\theta,x),
\;\;\;\; v(t,x)\rightarrow v(\theta,x)\frac{d\theta}{dt}.
\end{equation}
(Notice the different scaling for $v$, enforcing
$v(t,x)dt\rightarrow v(\theta,x)d\theta$.)
This leads, after short calculations, to
$$
\partial_\theta \rho+\nabla\cdot(\rho v)
=0,\;\;\;
\rho v+
2\theta[\partial_\theta (\rho v)+\nabla\cdot(\rho v\otimes v)]=
-\nabla  p(\rho).
$$
In the regime $\theta>>1$, we get the asymptotic model of ``pressureless" gas dynamics
$$
\partial_\theta \rho+\nabla\cdot(\rho v)
=0,\;\;\;
\partial_\theta (\rho v)+\nabla\cdot(\rho v\otimes v)=0,
$$
while, as ${\theta<<1}$, we recover
the Darcy law and the porous medium equation
$$
\rho v=-\nabla  p(\rho),\;\;\;
\partial_\theta \rho=\bigtriangleup (p(\rho)),
$$
and, in the isothermal case $p=\rho$, the heat equation
$$
\partial_\theta \rho=\bigtriangleup \rho.
$$

\subsection{From string motion to curve-shortening}
Let us now move to a model at the interface of geometry and high energy physics.
We consider a surface
$$
(t,s)\in\Omega \subset{\mathbb{R}}^2\rightarrow (t,X(t,s))\in
{\mathbb{R}}\times{\mathbb{R}}^d,
$$
parameterized by
a sufficiently smooth (at least Lipschitz continuous)
function $X$ over a bounded open space-time cylinder $\Omega$.
According to classical string theory (see \cite{Po}, for instance),
this surface is a relativistic string if and only if
$X$ is a critical point, with respect
to all smooth perturbations, compactly supported in $\Omega$,
of the ``Nambu-Goto Action" defined by
$$
\int_\Omega \sqrt{\partial_s X^2(1-\partial_t X^2)
+(\partial_t X\cdot\partial_s X)^2}\;\;
dtds
$$
which is nothing but the area of the surface,
in the space ${\mathbb{R}}\times{\mathbb{R}}^d$,
with respect to the Minkowski metric $(-1,+1,\cdot\cdot\cdot,+1)$.
It is customary to regularize this setting by
viewing
$$
(t,s)\in\Omega \subset{\mathbb{R}}^2
\rightarrow (t,s,X(t,s))\in
{\mathbb{R}^2}\times{\mathbb{R}}^d,
$$
as a graph in the enlarged space ${\mathbb{R}}^2\times{\mathbb{R}}^d$
and considering its area
in the enlarged Minkowski space ${\mathbb{R}}^2\times{\mathbb{R}}^d$,
with (rescaled) Minkowski metric $(-1,+\epsilon^2,+1,\cdot\cdot\cdot,+1)$:
$$
\int_\Omega \sqrt{(\epsilon^2+\partial_s X^2)(1-\partial_t X^2)
+(\partial_t X\cdot\partial_s X)^2}\;\;
dtds.
$$
[Of course, we recover the previous setting just as the special (and degenerate) case
$\epsilon=0$.]
The variational principle implies that $X$ is a solution to
the following
first order partial differential system (of hyperbolic type as $\epsilon>0$):
\begin{equation}
\label{string}
\partial_t(F\partial_t X-G\partial_s X)
-\partial_s(G\partial_t X+H\partial_s X)=0,
\end{equation}
where
$$
F=\frac{\epsilon^2+\partial_s X^2}{S},\;\;\;
G=\frac{\partial_t X\cdot\partial_s X}{S},\;\;\;
H=\frac{1-\partial_t X^2}{S},
$$
$$
S=\sqrt{(\epsilon^2+\partial_s X^2)(1-\partial_t X^2)
+(\partial_t X\cdot\partial_s X)^2}.
$$
After performing the quadratic change of time $\theta=t^2/2$, as we did in the previous
subsections, while keeping
only the zeroth order terms with respect to $\theta$, we easily obtain, as asymptotic equation
the nonlinear equation of parabolic type:
\begin{equation}
\label{curve shortening graph}
((\epsilon^2+\partial_s X^2)\mathbb{I}-\partial_s X\otimes\partial_s X)\partial_\theta X
=\sqrt{\epsilon^2+\partial_s X^2}\;\partial_s(\frac{\partial_s X}{\sqrt{\epsilon^2+\partial_s X^2}}).
\end{equation}
(Notice that $(\epsilon^2+\partial_s X^2)\mathbb{I}-\partial_s X\otimes\partial_s X$
is an invertible symmetric matrix with eigenvalues larger or equal to $\epsilon^2$.)
In the limit case $\epsilon=0$ , we get
\begin{equation}
\label{curve shortening gauge}
(\partial_s X^2\mathbb{I}-\partial_s X\otimes\partial_s X)\partial_\theta X
=|\partial_s X|\partial_s(\frac{\partial_s X}{|\partial_s X|})
\end{equation}
which becomes an ambiguous evolution equation, since it leaves
$\partial_\theta X\cdot\partial_s X$ undetermined.
[As a matter of fact, this geometric equation is not modified by any smooth
time-independent change of parameterization of the curve $s\rightarrow \sigma(s)$.]
\\
However,
we may solve instead the simpler equation
\begin{equation}
\label{curve shortening}
\partial_\theta X
=\frac{1}{|\partial_s X|}\partial_s(\frac{\partial_s X}{|\partial_s X|}).
\end{equation}
Indeed, this is a consistent way of solving (\ref{curve shortening gauge}) since
$$
\partial_\theta X\cdot \partial_s X=
\frac{\partial_s X}{|\partial_s X|}
\cdot\partial_s(\frac{\partial_s X}{|\partial_s X|})
=\partial_s(\frac{\partial_s X\cdot\partial_s X}{2|\partial_s X|^2})=0.
$$
Finally, by doing so, we have just recovered the familiar model of ``curve-shortening" in the Euclidean space $\mathbb{R}^d$ (see \cite{De}, for instance).

\subsection{The Eulerian form of the curve-shortening flow}

The string equation  (\ref{string}) admits a useful  ``Eulerian" version
\begin{equation}\label{eq:string-1}
\partial_{t}B + \nabla\cdot\left(\frac{B\otimes P-P\otimes B}{\rho}\right)=0,
\end{equation}
\begin{equation}\label{eq:string-2}
\nabla\cdot B=0,\;\;\;
\rho=\sqrt{B^2+P^2},
\end{equation}
\begin{equation}\label{eq:string-3}
\partial_{t}P + \nabla\cdot\left(\frac{P\otimes P}{\rho}\right)
= \nabla\cdot \left(\frac{B\otimes B}{\rho}\right)
\end{equation}
(which reads, in coordinates,
$$
\partial_t B^i+\partial_j(\rho^{-1}(B^i P^j-B^j P^i))=0,\:\:\:\partial_iB^i=0,
$$
$$
\rho=\sqrt{B_i B^i+P_i P^i},\;\;\;
\partial_tP^i+\partial_j(\rho^{-1}(P^i P^j-B^j B^i))=0).
$$
[As a matter of fact, defining
$$
B(t,x)=\int_\mathbb{R/Z}\delta(x-X(t,s))\partial_s X(t,s)ds,
$$
(which automatically satisfies $\nabla\cdot B=0$), assuming $X$ to be smooth, not self-intersecting,
with $\partial_s X$ never vanishing,
we get, after elementary calculations (similar to the ones done for the curve-shortening flow in Appendix 2, below),
that $B$ solves equations (\ref{eq:string-1},\ref{eq:string-2},\ref{eq:string-3})  together with
$$
P(t,x)
=\int_\mathbb{R/Z}\delta(x-X(t,s))\frac{
\left(\partial_s X^2\mathbb{I}-\partial_s X\otimes\partial_s X\right)\partial_t X(t,s)
}{
\sqrt{\partial_s X^2(1-\partial_t X^2)
+(\partial_t X\cdot\partial_s X)^2}}
ds.]
$$

Importantly enough, this system admits an extra conservation law:
\begin{equation}\label{eq:string-4}
\partial_t\rho+\nabla\cdot P=\nabla\cdot\left(\frac{(P\cdot B)B}{\rho^2}\right),\;\;\;\rho=\sqrt{B^2+P^2},
\end{equation}
which describes the local conservation of energy.
[This is easy to check. Indeed, using coordinates, we find
$$
\partial_t\rho=\frac{B_i\partial_t B^i+P_i\partial_t P^i}{\rho}
=\frac{B_i}{\rho}\partial_j\left(\frac{B^jP^i-B^i P^j}{\rho}\right)+\frac{P_i}{\rho}\partial_j\left(\frac{B^jB^i-P^i P^j}{\rho}\right)
$$
and notice that the second and fourth terms of the right-hand side combine as:
$$
-P^j\partial_j\left(\frac{P^2+B^2}{2\rho^2}\right)-\frac{P^2+B^2}{\rho^2}\partial_jP^j=-\partial_jP^j
$$
(since $\rho^2=B^2+P^2$), while the first and third terms give:
$$
\partial_j\left(\frac{P^i B^j B_i}{\rho^2}\right)
$$
(using $\nabla\cdot B=0$), which leads to the ``entropy conservation law" (\ref{eq:string-4}).]
\\
\\
Let us now perform the quadratic change of time:
$$
t\rightarrow \theta=\frac{t^2}{2},\;\;\;
B\rightarrow B,\;\;\;
P\rightarrow \frac{d\theta}{dt}P,\;\;\;
$$
which leads, as $\theta<<1$, to the asymptotic system
\begin{equation}\label{eq:short-1}
\partial_{\theta}B + \nabla\cdot\left(\frac{B\otimes P-P\otimes B}{\rho}\right)=0,
\end{equation}
\begin{equation}\label{eq:short-2}
\nabla\cdot B=0,\;\;\;\rho=|B|,
\end{equation}
\begin{equation}\label{eq:short-3}
P= \nabla\cdot \left(\frac{B\otimes B}{\rho}\right).
\end{equation}
Notice that this implies $B\cdot P=0$, since, in coordinates,
$$
B\cdot P=B_i\partial_j\left(\frac{B^i B^j}{\rho}\right)=\rho\frac{B_i}{\rho}B^j\partial_j\left(\frac{B^i}{\rho}\right)
=\rho B^j\partial_j\left(\frac{B^2}{2\rho^2}\right)=0
$$
(using $\nabla\cdot B=0$ and $|B|=\rho$).
We also get the extra equation, derived from (\ref{eq:string-4}),
\begin{equation}\label{eq:short-4}
\partial_\theta\rho+\frac{P^2}{\rho}+\nabla\cdot P=0,\;\;\;\rho=|B|.
\end{equation}
[Indeed, $\sqrt{B^2+2\theta P^2}=|B|+ \frac{\theta P^2}{|B|}+O(\theta^2)$,
which leads to
$$\partial_\theta \sqrt{B^2+2\theta P^2}=\partial_\theta |B|+ \frac{P^2}{|B|}+O(\theta).$$
We also used $B\cdot P=0$.]
Notice that this equation is no longer in conservation form, due to the emergence of the
dissipation term $P^2/\rho$ after the quadratic change of time  $t\rightarrow \theta$.
Equations (\ref{eq:short-1},\ref{eq:short-4}) also provide  a priori bounds for
any smooth solutions $B$ and $P$ on the flat torus $\mathbb{T}^d$
(that we have already taken into account in Definition \ref{def:dddf1}):
$$
\int_{\mathbb{T}^d}
 |B(\theta)|\le
\int_{\mathbb{T}^d} |B(0)|,\;\forall \theta\in [0,T],
\;\;\;\int_0^T\int_{\mathbb{T}^d} \frac{P^2}{|B|}\le \int_{\mathbb{T}^d}|B(0)|,
$$
$$
\int_{0}^{T}
 \int_{\mathbb{T}^d} |P|
 \le  \sqrt{
 \int_{0}^{T}
 \int_{\mathbb{T}^d} \frac{P^2}{|B|}}
 \sqrt{
\;\;  \int_{0}^{T}
 \int_{\mathbb{T}^d} |B|
 }
 \le \sqrt{T}\int_{\mathbb{T}^d}|B(0)|,
$$
and, for $0\le\theta_0\le\theta_1\le T$, and any smooth vector field $\varphi=\varphi(x)$,
$$
\left(\int_{\mathbb{T}^d} \left(B^i(\theta_1)-B^i(\theta_0)\right)\varphi_i \right)^2= \left(\int_{\theta_0}^{\theta_1}
\int_{\mathbb{T}^d} (\partial_j \varphi_i-\partial_i \varphi_j) \frac{P^j B^i}{|B|}\right)^2
$$
$$
 \le {\rm{Lip}}(\varphi)^2(\theta_1-\theta_0)
 \int_{\theta_0}^{\theta_1}
 \left(\int_{\mathbb{T}^d} |P|\right)^2
$$
(where ${\rm{Lip}}(\varphi)$ denotes the Lipschitz constant of $\varphi$),
$$
 \le {\rm{Lip}}(\varphi)^2(\theta_1-\theta_0)
 \int_{\theta_0}^{\theta_1}
 \left(\int_{\mathbb{T}^d} \frac{P^2}{|B|}\int_{\mathbb{T}^d}|B|\right)
$$
$$
 \le {\rm{Lip}}(\varphi)^2(\theta_1-\theta_0)
\left(\int_{\mathbb{T}^d} |B(0)|\right)^2,
$$
which shows that $B$ is bounded in $C^{1/2}([0,T],(C^1(\mathbb{T}^d))'_{w^*})$
by a constant depending only on $\int_{\mathbb{T}^d}|B(0)|$ and $T$.
\\
\\
Equations (\ref{eq:short-1},\ref{eq:short-2},\ref{eq:short-3})
can also be written in non-conservative form
in terms of
$$
b=\frac{B}{\rho},\;\;\;v=\frac{P}{\rho},\;\;\;\rho=|B|.
$$
We already have $b^2=1$ and $b\cdot v=0$.
Using coordinates,
we first get from (\ref{eq:short-3})
$$
\rho v^i=P^i=\partial_j\left(\frac{B^i B^j}{\rho}\right)=\partial_j(\rho b^i b^j)=\rho b^j\partial_j b^i
$$
(since $\partial_j(\rho b^j)=\partial_jB^j=0$).
Next, (\ref{eq:short-1}) becomes
$$
\rho(\partial_\theta b^i+ v^j\partial_j b^i-b^j\partial_j v^i)=-b^i(\partial_\theta\rho+\partial_j(\rho v^j))
=b^i \rho v^2
$$
(thanks to (\ref{eq:short-4})).
So we have obtained
\begin{equation}
\partial_\theta b+(v\cdot\nabla)b=(b\cdot\nabla)v+b v^2,
\;\;\;v=(b\cdot\nabla)b,
\end{equation}
(which is consistent with $b^2=1$ and $b\cdot v=0$ as can be easily checked).
Notice that (\ref{eq:short-4}) can be written according to the non-conservative variables as
\begin{equation}\label{eq:short-5}
\partial_\theta\rho+\nabla\cdot(\rho v)=-\rho v^2,
\end{equation}
which is a linear equation in $\rho$.
\section{Analysis of the Eulerian curve-shortening flow}
\label{main}
\subsection{Relative entropy for the Eulerian equations for strings}

We start from the ``Eulerian'' version
(\ref{eq:string-1},\ref{eq:string-2},\ref{eq:string-3}) of the string equation
(\ref{string}).
This system belongs to the class of systems of conservation laws:
$$
\partial_t V+\nabla\cdot \mathcal{F}(V)=0,
$$
where $\mathcal{F}$ is a given function and $V$ is a vector-valued function (for us $V=(B,P)$).
If such a system admits an $extra$ conservation law
$$
\partial_t \mathcal{E}(V)+\nabla\cdot \mathcal{G}(V)=0,
$$
for a pair of functions $\mathcal{E,G}$, with $\mathcal{E}$ strictly convex, then the system is automatically
``hyperbolic" (i.e. well posed, at least for short time), under minor additional conditions \cite{Da}, and $\mathcal{E}$
is often called an ``entropy" for the system (although it should be called ``energy" for a large class of applications).
The system for strings (\ref{eq:string-1},\ref{eq:string-2},\ref{eq:string-3}) admits such an extra conservation law, namely
(\ref{eq:string-4}), with $V=(B,P)\rightarrow \mathcal{E}(V)=\sqrt{B^2+P^2}$ as entropy function.
[Notice, however, that the entropy $\rho$ is not strictly convex.]

Hyperbolic systems of conservation laws with convex entropy enjoy
a  ``weak-strong uniqueness" principle \cite{Da},
based on the concept of ``relative entropy" (also called  ``modulated energy" or ``Bregman divergence", depending
of the frameworks):
$$
\eta(V,V^*)=\mathcal{E}(V)-\mathcal{E}(V^*)-\nabla\mathcal{E}(V^*)\cdot(V-V^*),
$$
which is just the discrepancy between $\mathcal{E}$ at point $V$ and its linear approximation about a given point $V^*$.
(Observe that, as $\mathcal{E}$ is a convex function with Hessian bounded away from zero and infinity, the relative entropy behaves
as $|V-V^*|^2$.) Notice that the relative entropy is as convex as the entropy as a function of $V$ ($V^*$ being kept fixed)
since it differs just by an affine term.
\\
\\
In the case of system (\ref{eq:string-1},\ref{eq:string-2},\ref{eq:string-3}), the relative entropy density is defined, for $(B,P)\in\mathbb{R}^{d}\times \mathbb{R}^{d}$
and
$(b^*,v^*)\in\mathbb{R}^{d}\times \mathbb{R}^{d}$, by
$$
\eta=\frac{{|B-\rho b^*|}^2+{|P-\rho v^*|}^2}{2\rho}=
\rho\frac{1+{b^*}^2+{v^*}^2}{2}-B\cdot b^*-P\cdot v^*,\;\;\;\rho=\sqrt{B^2+P^2},
$$
which is convex in $(B,P)$. (Notice that, at this stage, we do not assume ${b^*}^2+{v^*}^2=1$, which would be natural
to define the relative entropy but would lead to contradictions after performing the quadratic change of time as will be done
in the next subsection.)
\\
Let us now consider a smooth, $\mathbb{Z}^d-$ periodic in space, solution $(B,P)(t,x)$ of equations
(\ref{eq:string-1},\ref{eq:string-2},\ref{eq:string-3})
and monitor the evolution, on a fixed time interval $[0,T]$, of
the integral of $\eta$ over ${\mathbb{(R/Z)}^d}$,
for some smooth trial functions:
$$
(t,x)\in[0,T]\times\mathbb{(R/Z)}^d
\rightarrow (b^*(t,x)\in\mathbb{R}^d,\;v^*(t,x)\in\mathbb{R}^d).
$$
After tedious and elementary calculations, we find
\begin{equation}\label{eq:in-1}
\begin{array}{r@{}l}
\displaystyle{\quad \frac{d}{dt}\int \eta\;=} & \displaystyle{
\; \int\frac{1}{2\rho}
(B_i-\rho b^*_i)(B_j-\rho b^*_j)(\partial_jv^*_{i}+\partial_iv^*_{j})
} \\
& \displaystyle{-\int\frac{1}{2\rho}
(P_i-\rho v^*_i)(P_j-\rho v^*_j)(\partial_jv^*_{i}+\partial_iv^*_{j})} \\
& \displaystyle{- \int\frac{1}{\rho}
(B_i-\rho b^*_i)(P_j-\rho v^*_j)
(\partial_jb^*_{i}-\partial_ib^*_{j}) }\\
& \displaystyle{+ \int \eta\cdot \mathrm{L}_1+ \int B\cdot \mathrm{L}_2+ \int P\cdot \mathrm{L}_3  -\int\frac{(P\cdot B)B}{\rho^2}\cdot\nabla\left(\frac{{b^*}^2+{v^*}^2}{2}\right) }\\
\end{array}
\end{equation}
where
$$
\mathrm{L}_1=\frac{\zeta^*}{1+{b^*}^2+{v^*}^2},$$
$$
\zeta^*=D_t^*({b^*}^2+{v^*}^2)-2b^*\cdot \nabla(b^*\cdot v^*),  
\;\;\;D^*_t=(\partial_t+v^*\cdot\nabla),
$$
$$
\mathrm{L}_2=-D^*_t b^*+(b^*\cdot\nabla)v^*+\nabla(b^*\cdot v^*)+b^*{\rm{L}}_1
$$
$$
\;\;\;\mathrm{L}_3=-D^*_t v^*+(b^*\cdot\nabla)b^*+v^*{\rm{L}}_1 
$$
\subsection{Relative entropy and quadratic change of time}
After the quadratic change of time,
$$
t\rightarrow \theta=\frac{t^2}{2},\;\;\;
\rho\rightarrow \rho,\;\;\;
B\rightarrow B,\;\;\;b^*\rightarrow b^*,\;\;\;
P\rightarrow \theta'(t)
P,\;\;\;
v^*\rightarrow \theta'(t)
v^*,
$$
we get
$\theta'(t)^2=2\theta,\;\;\;\theta''(t)=1$,
$$
\begin{array}{r@{}l}
\displaystyle{\quad \frac{d}{d\theta}\int \eta\;= }& \displaystyle{\;\int\frac{1}{2\rho}
(B_i-\rho b^*_i)(B_j-\rho b^*_j)(\partial_jv^*_{i}+\partial_iv^*_{j})
} \\
& \displaystyle{
-\int\frac{\theta}{\rho}
(P_i-\rho v^*_i)(P_j-\rho v^*_j)(\partial_jv^*_{i}+\partial_iv^*_{j})} \\
& \displaystyle{- \int\frac{1}{\rho}
(B_i-\rho b^*_i)(P_j-\rho v^*_j)
(\partial_jb^*_{i}-\partial_ib^*_{j}) }\\
& \displaystyle{ + \int \eta\cdot \mathrm{L}_1+ \int B\cdot \mathrm{L}_2+ \int P\cdot \mathrm{L}_3 -\int\frac{(P\cdot B)B}{\rho^2}\cdot\nabla\left(\frac{{b^*}^2+2\theta{v^*}^2}{2}\right), }\\
\end{array}
$$
where
$$
\mathrm{L}_1=\left(1+{b^*}^2+2\theta{v^*}^2\right)^{-1}\zeta^*,$$
$$
\zeta^*=D_\theta({b^*}^2+2\theta{v^*}^2)-2b^*\cdot \nabla(b^*\cdot v^*)
,\;\;\;
D^*_\theta=(\partial_\theta+v^*\cdot\nabla),
$$
$$
\mathrm{L}_2=-D^*_\theta b^*+(b^*\cdot\nabla)v^*+\nabla(b^*\cdot v^*)+b^*{\rm{L}}_1
$$
$$
\;\;\;\mathrm{L}_3=-v^*-2\theta D^*_\theta v^*+(b^*\cdot\nabla)b^*
+2\theta v^*{\rm{L}}_1.
$$
Now, in order to address the Eulerian curve-shortening system,
we want to drop the terms of order $O(\theta)$ and limit ourself to the case when ${b^*}^2=1$.
However, we have to be very careful about all terms involving $\partial_\theta$.
This happens first in the definition of $\zeta^*$, because of the term
$$
D_\theta({b^*}^2+2\theta{v^*}^2)=D_\theta({b^*}^2)+2{v^*}^2+O(\theta).
$$
So, in the limit $\theta=0$, with ${b^*}^2=1$,
we find
$$
\zeta^*=2{v^*}^2-2b^*\cdot \nabla(b^*\cdot v^*)
$$
and, therefore,
$$
\mathrm{L}_1={v^*}^2-b^*\cdot \nabla(b^*\cdot v^*).
$$
Similarly, we have to take care of
$$
\frac{d}{d\theta}\int \eta.
$$
where
$$
\eta=\frac{1+{b^*}^2+2\theta{v^*}^2}{2}\rho-B\cdot b^*-2\theta P\cdot v^*
$$
and
$$
\rho=\sqrt{B^2+2\theta P^2}
=|B|+\theta\frac{P^2}{|B|}+O(\theta^2)
$$
We get
$$
\partial_\theta \rho
=\partial_\theta|B|+\frac{P^2}{|B|}+O(\theta),
$$
$$
\partial_\theta \eta
={v^*}^2\rho +\partial_\theta(\frac{1+{b^*}^2}{2}\rho)-
\partial_\theta(B\cdot b^*)-2P\cdot v^*+O(\theta)
$$
This leads, as ${b^*}^2=1$, to
$$
\rho=|B|+O(\theta),\;\;\;\eta=|B|-B\cdot b^*+O(\theta),
$$
$$
\partial_\theta \eta
={v^*}^2\rho +\partial_\theta|B|+\frac{P^2}{|B|}
-\partial_\theta(B\cdot b^*)-2P\cdot v^*+O(\theta)
$$
$$
=\partial_\theta(|B|-B\cdot b^*)
+\frac{(P-|B|v^*)^2}{|B|}
+O(\theta)
=\partial_\theta(\rho-B\cdot b^*)
+\frac{(P-\rho v^*)^2}{\rho}
+O(\theta).
$$
Finally, after dropping the terms of order $O(\theta)$ and limiting ourself to the case when ${b^*}^2=1$,
we have found
$$
\mathrm{L}_1={v^*}^2-b^*\cdot \nabla(b^*\cdot v^*),
$$
$$
\mathrm{L}_2=-D^*_\theta b^*+(b^*\cdot\nabla)v^*+\nabla(b^*\cdot v^*)+b^*{\rm{L}}_1
$$
$$
=-D^*_\theta b^*+(b^*\cdot\nabla)v^*+\nabla(b^*\cdot v^*)+b^*({v^*}^2-b^*\cdot \nabla(b^*\cdot v^*)),
$$
$$
D^*_\theta=(\partial_\theta+v^*\cdot\nabla),\;\;\;\mathrm{L}_3=-v^*+(b^*\cdot\nabla)b^*,
$$
and, for all smooth trial field $b^*$ such that ${b^*}^2=1$,
\begin{equation}\label{eq:in-1}
\begin{array}{r@{}l}
\displaystyle{\quad \frac{d}{d\theta}\int \eta
+\int\frac{(P-\rho v^*)^2}{\rho}\;=} & \displaystyle{\;
\int\frac{1}{2\rho}
(B_i-\rho b^*_i)(B_j-\rho b^*_j)(\partial_jv^*_{i}+\partial_iv^*_{j})
} \\
& \displaystyle{- \int\frac{1}{\rho}
(B_i-\rho b^*_i)(P_j-\rho v^*_j)
(\partial_jb^*_{i}-\partial_ib^*_{j}) }\\
& \displaystyle{
+ \int \eta\; \mathrm{L}_1+ \int B\cdot \mathrm{L}_2+ \int P\cdot \mathrm{L}_3 },\\
\end{array}
\end{equation}
where
$$
\rho=|B|,\;\;\;\eta=\rho-B\cdot b^*=|B|-B\cdot b^*=\frac{(B-|B|b^*)^2}{2|B|}=\frac{(B-\rho b^*)^2}{2\rho}.
$$
Using the Cauchy-Schwarz inequality, we can find a constant $c^*$ depending only on $b^*$ and $v^*$ such that
$$
\begin{array}{r@{}l}
& \displaystyle{\quad \frac{d}{d\theta}\int \eta
+\int\frac{|P-\rho v^*|^2}{2\rho}
\le c^* \int \eta\;dx
+ \int B\cdot \mathrm{L}_2+ \int P\cdot \mathrm{L}_3 }.\\
\end{array}
$$
This implies, for any constant $r\ge c^*$,
$$
(-r+\frac{d}{d\theta})\int \eta
+\int\frac{|P-\rho v^*|^2}{2\rho}+ (r-c^*) \int \eta\;dx
\le \int B\cdot \mathrm{L}_2+ \int P\cdot \mathrm{L}_3
$$
and, after multiplying this inequality by $e^{-r\theta}$ and integrating in time
 $\sigma\in[0,\theta]$,
\begin{equation}\label{dissipative-0}
 e^{-r\theta
 }\int \eta(\theta)
+\int_0^\theta e^{-r \sigma}
\left((r-c^*)\int \eta+\int \frac{|P-\rho v^*|^2}{2\rho}-R\right)
(\sigma)d\sigma\le \int \eta(0),
\end{equation}
where
$$
R=\int B\cdot \mathrm{L}_2+ \int P\cdot \mathrm{L}_3,
$$
$$
\mathrm{L}_2=-D^*_\theta b^*+
(b^*\cdot\nabla)v^* +b^*{v^*}^2 -b^*(b^*\cdot\nabla)(b^*\cdot v^*),
\;\;\;\mathrm{L}_3=-v^*+(b^*\cdot\nabla)b^*.
$$
We can write
\begin{equation*}
\begin{array}{r@{}l}
\displaystyle{ \frac{|P-\rho v^*|^2}{2\rho}\; }& \displaystyle{=\; \sup_{A}\;(P-\rho v^*)\cdot A -\rho\frac{A^2}{2}}\\
&\displaystyle{=\; \sup_{A} P\cdot A -(\eta + B\cdot b^*)\left(v^*\cdot A+\frac{A^2}{2}\right)}
\end{array}
\end{equation*}
(since $\eta=\rho-B\cdot b^*$)
and substitute for inequality (\ref{dissipative-0})
the family of inequalities
\begin{multline}
\label{dissipative-1-bis}
e^{-r\theta}
\int \eta(\theta)
+\int_0^\theta e^{-r\sigma}
\left[\int P\cdot (A-L_3)+\left(r-c^*-\frac{A\cdot(A+2v^*)}{2}\right)\eta\right.\\
\left.-B\cdot \left(L_2+b^*\frac{A\cdot(A+2v^*)}{2}\right)\right]
(\sigma)d\sigma\le \int \eta(0).
\end{multline}
Observe that these inequalities are convex in $(B,P)$ as long as $r$ is chosen so that
$$
r\ge c^*+\sup_{\theta,x}\frac{A\cdot(A+2v^*)}{2}.
$$

However, this creates a problem, since $r$ must depend on $A$.
This is why we input
a cut-off parameter
$\lambda>0$ and assume that the trial functions $A$ are chosen with $|A(\theta,x)|\le \lambda$.
By doing this, the advantage is that we maintain the convexity of inequality as long as $r$ is chosen
big enough only as a function of $b^*,v^*$ and $\lambda$,
namely:
$$
r\ge c^*+\frac{\lambda^2}{2}+\lambda\|v^*\|_{\infty}.
$$
The price to pay is that we cannot fully recover
$$
\frac{|P-\rho v^*|^2}{2\rho}
$$
by taking the supremum over all $A$ such that $|A|\le \lambda$, but only the $\lambda-$approximation
$
K_\lambda(\rho,P-\rho v^*),
$
where
$$
K_\lambda(\rho,Z)
=\sup_{|A|\le \lambda} Z\cdot A-\rho\frac{A^2}{2}
=\frac{|Z|^2}{2\rho}-\frac{(|Z|-\lambda\rho)_+^2}{2\rho}\ge 0.
$$
Observe that, by doing so, we keep a good control of the distance between $P$ and $\rho v^*$, since
(as can be easily checked)
\begin{equation}
\label{K-inequality}
K_\lambda(\rho,P-\rho v^*)\ge \min\left(\frac{(P-\rho v^*)^2}{2\rho},\frac{\lambda|P-\rho v^*|}{2}\right).
\end{equation}
So, the supremum of inequalities (\ref{dissipative-1-bis})
over all trial functions $A$ such that $|A|\le \lambda$,
is $equivalent$ to
\begin{equation}
\label{dissipative-3-bis}
e^{-r\theta}
\int \eta(\theta)
+\int_0^\theta e^{-r\sigma}
(\int K_\lambda(\rho,P-\rho v^*)+(r-c^*)\eta
-R)
(\sigma)d\sigma\le \int \eta(0).
\end{equation}
Now let us consider $(B,P)$ not only as functions but also as vector-valued Borel measures,
for which (\ref{dissipative-1-bis}) is still well-defined. The $\lambda-$approximation $K_\lambda(\rho,P-\rho v^*)$ can be interpreted as a function of measures \cite{DeTe} and (\ref{dissipative-3-bis}) is equivalent to (\ref{dissipative-1-bis}) in the sense that,
$$\int_0^\theta e^{-r\sigma}\int K_\lambda(\rho,P-\rho v^*)=\sup_{A\in C^{0}\atop\|A\|_{\infty}\le\lambda} \int_0^\theta e^{-r\sigma}\int (P-\rho v^*)\cdot A -\rho\frac{A^2}{2}.$$
Notice that, due to the convexity of $K_\lambda$,
we have
\begin{equation}\label{K-in}
 \int_0^\theta e^{-r\sigma}\int K_\lambda(\rho,P-\rho v^*)(\sigma)d\sigma
 \geq e^{-r\theta} \int_0^\theta  K_\lambda\left(\int \rho(\sigma),\int |P-\rho v^*|(\sigma)\right)d\sigma.
\end{equation}
With these calculations, we have recovered
the concept of dissipative solutions as given in Definition \ref{def:dddf1}.
Then,
the proof of our main results becomes
straightforward.

\subsubsection*{Proof of Theorem \ref{theorem1}}

We just have to show that, for fixed initial conditions $B_0$, the set of dissipative solutions,
as defined by Definition \ref{def:dddf1}, if not empty, is convex and compact
for the weak-* topology of measures. The convexity of the set of solutions is almost free. It follows directly from the convexity of inequalities (\ref{dissipative-1}). Let's focus on the compactness. Our goal is to prove that, if $\{(B_n,P_n)\}_{n\in\mathbb{N}}$ is a sequence of dissipative solutions with initial data $B_0$, then up to a subsequence, it converges in the weak-* topology of measures to a dissipative solution $(B,P)$ with the same initial data. This follows from the inequalities (\ref{dissipative-1}) and suitable bounds that we assume for $B_n$ and $P_n$. To see this, let's first show that, $\sup_{\theta}\int|B_n(\theta)|$ is uniformly bounded. (Indeed, let's take $b^*=(1,0,\ldots,0)$, $v^*=A=0$ in (\ref{dissipative-1}). Then we have $\int |B_n(\theta)|-B^1_{n}(\theta)\leq C\int |B_0|$, $\forall\theta\in[0,T]$. Since $B_n$ is bounded in $C^{1/2}([0,T],(C^{1}(\mathbb{T}^d))'_{w^*})$, there exists a constant $C'$ such that for any $n,\theta$, $|\int B^1_{n}(\theta)-\int B_0^1|\le C'$. So we get a uniform upper bound of $\int |B_n(\theta)|$.) Therefore, for any $\theta\in[0,T]$, the set $\{B_n(\theta)\}_{n\in\mathbb{N}}$ is relatively compact for the weak-* topology of $C(\mathbb{T}^d
,\mathbb{R}^{d})'$. Next, we look at the map $[0,T]\rightarrow C(\mathbb{T}^d
,\mathbb{R}^{d})'_{w^*}$, $t\rightarrow B(t)$. This map is equicontinuous because of the assumption on ${B_n}$. Then, by Arzel\`{a}-Ascoli's theorem, there exists $B\in C([0,T],C(\mathbb{T}^d,\mathbb{R}^{d})'_{w^{*}})$, such that, up to a subsequence, $B_n(\theta)\rightharpoonup^{*}B(\theta)$, $\forall\theta\in[0,T]$. Now, since $\iint |P_n|$ is bounded, there exists $P\in C([0,T]\times \mathbb{T}^d,\mathbb{R}^{d})'$, such that $P_n\rightharpoonup^{*}P$. Then because inequalities (\ref{dissipative-1}) are stable under weak-* convergence, we can prove that the limit $(B,P)$ satisfies all the requirements in Definition \ref{def:dddf1}, therefore, it is also a dissipative solution with initial data $B_0$.

\subsubsection*{Proof of Theorem \ref{theorem2}}

Let $(b,v)$ be a smooth solution of the non-conservative
form of the curve-shortening flow (\ref{eq:shortnc}) with $b^2=1$, which directly implies
$b\cdot v=0$.
We have to show that any dissipative solution satisfies $B=|B| b$ and $P=|B| v,$ as soon as
$B(0)=|B(0)|b(0)$.
The proof is quite straightforward:
we already have $\eta(0)=0$ since $B(0)=|B(0)|b(0)$.
Next, we set $b^*=b$, $v^*=v$, $A=0$ and fix $\lambda>0$ in definition (\ref{dissipative-1}).
Since we have (\ref{eq:shortnc}) and $b\cdot v=0$, we get ${\rm{L}}_2={\rm{L}}_3=0$.
Since $\eta\ge0$, the inequality \ref{dissipative-1} directly implies $\eta=0$, $\forall\theta\in[0,T]$,
and, therefore $B=|B|b$. Now, let's go back to the inequality \ref{dissipative-1} which is already simplified since $\eta=0$. By taking the supremum over all $A$ such that $\|A\|_{\infty}\le \lambda$, we get
$$\int_0^T e^{-r\sigma}\int K_\lambda(\rho,P-|B| v)\le 0.$$
Using (\ref{K-in}) we deduce
$$
\int_0^T  K_\lambda\left(\int \rho(\sigma),\int \big|P-|B|v\big|(\sigma)\right)d\sigma=0,
$$
 and, therefore, $P=|B|v$
(because of (\ref{K-inequality})), which completes the proof.

\section{
Appendix  1: modulated energy and dissipative solutions for ordinary dynamical equations
}
In this appendix, we
explain, in the very elementary case of our dynamical system, the concepts of ``modulated energy" (also called ``relative entropy") and ``dissipative formulation",
which will later be used and extended to the dissipative setting.
\\
Here, we crucially assume that the potential $\varphi$ is convex and, in order to keep
the presentation simple, we assume that the spectrum of the symmetric matrix $D^2\varphi(x)$ is
uniformly contained in some fixed interval $[r,r^{-1}]$ for some constant $r\in (0,1/2)$. We further
assume that the third derivatives of $\varphi$ are bounded.
The total energy of a curve $t\rightarrow X(t)$ is defined by
$$
\frac{1}{2}|X'(t)|^2+\varphi(X(t))\;\;\;({\rm where}\;\;X'(t)=\frac{dX}{dt}\;)
$$
and is a constant as $X$ is solution to the dynamical system
$$
X''(t)=-\nabla\varphi(X(t)).
$$
Given a smooth curve $t\rightarrow Y(t)$, we define the ``modulated energy" (or ``relative entropy") of $X$ at time $t$ with respect to $Y$
by expanding the energy about $Y$ at $X$:
$$
\eta[t,X,Y]=\frac{1}{2}|X'(t)-Y'(t)|^2+\varphi(X(t))-\varphi(Y(t))-\nabla\varphi(Y(t))\cdot(X(t)-Y(t)).
$$
Because of the assumption we made on $\varphi$, $\eta$ is a perfect substitute for the
squared distance between $(X,X')$ and $(Y,Y')$:
$$
r\le \frac{2\eta[t,X,Y]}{|X-Y|^2+|X'-Y'|^2}\le r^{-1}.
$$
We get
$$
\frac{d}{dt}\eta[t,X,Y]=(X'-Y')\cdot(
X''-Y'')
+\nabla\varphi(X)\cdot X'
-\nabla\varphi(Y)\cdot Y'
$$
$$
-\nabla\varphi(Y)\cdot(X'-Y')
-Y'\cdot D^2\varphi(Y)\cdot (X-Y)
$$
$$
=(X'-Y')\cdot(X''+\nabla\varphi(X)-Y''-\nabla\varphi(Y))+Y'\cdot(\nabla\varphi(X)-\nabla\varphi(Y)-D^2\varphi(Y)\cdot (X-Y)).
$$
We first observe that
$$
|\nabla\varphi(X)-\nabla\varphi(Y)-D^2\varphi(Y)\cdot (X-Y)|\le C|X-Y|^2\le C\eta[t,X,Y]
$$
where, from now on,
$C$ is a generic constant that depends only on $\varphi$ or $Y$.
So,
$$
\frac{d}{dt}\eta[t,X,Y]-
(X'-Y')\cdot(X''+\nabla\varphi(X)-Y''-\nabla\varphi(Y))
\le C\eta[t,X,Y]
$$
and then, after integration in time for $t\in [0,T]$, $T>0$ being an arbitrarily chosen fixed time,
\begin{equation}
\label{ineq}
\eta[T,X,Y]
-\int_0^T (X'(t)-Y'(t))\cdot (\omega_X(t)-\omega_Y(t))e^{(T-t)C}dt\le \eta[0,X,Y]e^{CT},
\end{equation}
where
$$
\omega_Z(t)=Z''(t)+\nabla\varphi(Z(t))
$$
\\
Let us exploit inequality (\ref{ineq}) in several different ways.
\\
First, we see that for a curve $X$ it is equivalent to be solution of the dynamical system, i.e. $\omega_X=0$
or to satisfy
\begin{equation}
\label{solution dissipative}
\eta[T,X,Y]
+\int_0^T (X'(t)-Y'(t))\cdot \omega_Y(t)e^{(T-t)C}dt\le \eta[0,X,Y]e^{CT},\;\;\;\forall T>0,
\end{equation}
for any smooth curve $Y$, where $C$ is a constant depending only on $Y$ (up to time $T$) and $\varphi$.
Indeed, by taking as $Y$ the unique solution of the dynamical system with initial conditions $Y(0)=X(0)$, $Y'(0)=X'(0)$
provided by the Cauchy-Lipschitz theorem on ODEs,
we get both $\omega_Y=0$ and $\eta[0,X,Y]=0$. Thus inequality (\ref{solution dissipative}) just says $\eta[T,X,Y]=0$ for
all $T>0$,
which means $X=Y$ and, therefore, $X$ is indeed a solution to the dynamical system.
Thus, we can take (\ref{solution dissipative}) as an alternative
notion of solution, that we call ``dissipative solution". This inequality has the advantage to be convex in $X$,
as the initial conditions $X(0),X'(0)$ are fixed, and therefore preserved under weak convergence of $(X,X')$.
\\
\\
Next, we use (\ref{solution dissipative}) to compare a solution $X$ of the dynamical system with zero initial velocity,
i.e. $X'(0)=0$,
to the solution $Z$ of the gradient flow equation
$$
Z'(\theta)+\nabla\varphi (Z(\theta))=0,
$$
with initial condition $Z(0)=X(0)$. Indeed, let us set $Y(t)=Z(\theta)$, $\theta=t^2/2$.
Then $Y'(t)=tZ'(\theta)$, $Y'(0)=0$, $Y(0)=Z(0)=X(0)$,  $Y''(t)=Z'(\theta)+t^2 Z''(\theta)$, which implies
$ \eta[0,X,Y]=0$ and $\omega_Y(t)=t^2 Z''(\theta)$. So, (\ref{solution dissipative}) gives
$$
\eta[T,X,Y]
+\int_0^T (X'(t)-Y'(t))\cdot t^2 Z''(t^2/2)e^{(T-t)C}dt\le 0,
$$
which implies (by Cauchy-Schwarz inequality and by definition of $\eta$)
$$
\eta[T,X,Y]
\le C\int_0^T (\eta[t,X,Y]+t^4)
dt,
$$
where $C$ is a generic constant depending only on $T$, $\varphi$ and $Z$. By Gronwall's lemma,
we conclude that $\eta[t,X,Y]\le Ct^5$ which implies, by definition of $\eta$,
$$
|X(t)-Z(t^2/2)|^2+|\frac{dX}{dt}(t)-t\frac{dZ}{d\theta}(t^2/2)|^2\le Ct^5,\;\;\;\forall t\in [0,T],
$$
as already claimed, at the beginning of this subsection. (See (\ref{error}).)
(Notice that the smallest expected error is $O(t^6)$ as shown by the example
$d=1$, $\varphi(x)=|x|^2/2$, for which $X(t)=X(0)\cos(t)$, while $Z(\theta)=X(0)\exp(-\theta)$.)
\section{Appendix 2: direct recovery of the Eulerian curve-shortening flow}

For the sake of completeness, let us check that system (\ref{eq:short-1},\ref{eq:short-2},\ref{eq:short-3})
indeed describes the curve-shortening
flow in $\mathbb{R}^d$, for a continuum of non intersecting curves.
Let us do the calculation in the case of a single smooth time-dependent loop,
$s\in\mathbb{R/Z}\rightarrow X(\theta,s)$, that we assume to be
non self-intersecting at every fixed time $\theta$, and such that
$\partial_s X$ never vanishes.
We introduce (as a distribution, or, if one prefers, as a ``$1-$current")
$$
B(\theta,x)=\int_\mathbb{R/Z}\delta(x-X(\theta,s))\partial_s X(\theta,s)ds,
$$
which automatically satisfies $\nabla\cdot B=0$.
Since $X$ is smooth, not self-intersecting, and $\partial_s X$ never vanishes,
by assumption,
we may find a smooth vector field $v(\theta,x)$
such that
$$
\partial_\theta X(\theta,s)=v(\theta,X(\theta,s))
$$
that we can interpret as the ``Eulerian velocity field" attached to the loop
evolution.
We also introduce the nonnegative field
$$
\rho(\theta,x)
=\int_\mathbb{R/Z}\delta(x-X(\theta,s))|\partial_s X(\theta,s)|ds
$$
which can also be interpreted as $|B(\theta,x)|$ since $X$ is supposed to be non self-intersecting.
We get (using indices $i,j,k\in\{1,\cdot\cdot\cdot,d\}$ with implicit summation on repeated indices)
$$
\partial_\theta B^i(\theta,x)=\int_\mathbb{R/Z}[
-(\partial_j\delta)(x-X(\theta,s))\partial_\theta X^j(\theta,s)\partial_s X^i(\theta,s)
$$
$$
+\delta(x-X(\theta,s))\partial^2_{s\theta} X^i(\theta,s)]ds
$$
(in distributional sense)
$$
=-\int_\mathbb{R/Z}
(\partial_j\delta)(x-X(\theta,s))
[\partial_\theta X^j(\theta,s)\partial_s X^i(\theta,s)
-\partial_\theta X^i(\theta,s)\partial_s X^j(\theta,s)]ds
$$
(after integration by part in $s\in \mathbb{R/Z}$ of the second term)
$$
=-\partial_j \int_\mathbb{R/Z}\delta(x-X(\theta,s))
[\partial_\theta X^j(\theta,s)\partial_s X^i(\theta,s)
-\partial_\theta X^i(\theta,s)\partial_s X^j(\theta,s)]ds.
$$
So
$$
\partial_\theta B(\theta,x)=-\nabla\cdot \int_\mathbb{R/Z}\delta(x-X(\theta,x))
(\partial_s X(\theta,s)\otimes \partial_\theta X(\theta,s)
-\partial_\theta X(\theta,s)\otimes \partial_s X(\theta,s))ds.
$$
Then we can write
$$
\partial_\theta B(\theta,x)+\nabla\cdot \int_\mathbb{R/Z}\delta(x-X)
(\partial_s X(\theta,s)\otimes v(\theta,X(\theta,s))
-v(\theta,X(\theta,s))\otimes \partial_s X(\theta,s))ds=0
$$
which exactly means, by definition of $B$,
\begin{equation}
\label{induction}
\partial_\theta B+\nabla\cdot(B\otimes v-v\otimes B)=0.
\end{equation}

Since $X$ is assumed to be non-intersecting,
by definition of $v$, we may write
$$
(|B|v)(\theta,x)=\int_\mathbb{R/Z}\delta(x-X(\theta,s))|\partial_s X(\theta,s)|v(\theta,X(\theta,s))ds
$$
So far, we have not used equation (\ref{curve shortening}), namely
$$
\partial_\theta X
=\frac{1}{|\partial_s X|}\partial_s(\frac{\partial_s X}{|\partial_s X|}),
$$
Let us do it now:
$$
(\rho v^i)(\theta,x)=(|B|v^i)(\theta,x)=\int_\mathbb{R/Z}\delta(x-X(\theta,s))\partial_s(\frac{\partial_s X^i}{|\partial_s X|})ds
$$
$$
=\int_\mathbb{R/Z}(\partial_j\delta)(x-X(\theta,s))\frac{\partial_s X^j\partial_s X^i}{|\partial_s X|}ds
$$
(after integrating by part in $s\in \mathbb{R/Z}$)
$$
=\partial_j\int_\mathbb{R/Z}\delta(x-X(\theta,s))\frac{\partial_s X^j\partial_s X^i}{|\partial_s X|}ds
$$
that we can interpret as
$$
\rho v=\nabla\cdot \frac{B\otimes B}{|B|}=\nabla\cdot \frac{B\otimes B}{\rho}.
$$
Finally we can write (\ref{induction}) as
$$
\partial_\theta B+\nabla\cdot(\frac{B\otimes P-P\otimes B}{\rho})=0,\;\;\;\nabla\cdot B=0,\;\;\;P=\nabla\cdot \frac{B\otimes B}{\rho},\;\;\;\rho=|B|,
$$
where $P$ stands for $\rho v$
and $(B,\rho,P)$ solves equations (\ref{eq:short-1},\ref{eq:short-2},\ref{eq:short-3}).
So far, our claim has been justified only in the case of a single loop.
We next argue that, due to its homogeneity of degree 1, equations (\ref{eq:short-1},\ref{eq:short-2},\ref{eq:short-3}),
in spite of their nonlinearity, enjoy a nice superposition
principle, in the sense that we may still get a solution by superposing several
smooth curves subject to curve-shortening as long as they do not intersect and we may even build
smooth solutions by using a continuum of such curves.
This concludes the proof of our claim
that equations (\ref{eq:short-1},\ref{eq:short-2},\ref{eq:short-3}),
are the ``Eulerian formulation" of the curve-shortening flow.
Notice that similar calculations can also be performed to justify
the Eulerian version
(\ref{eq:string-1},\ref{eq:string-2},\ref{eq:string-3}) of the string equation
(\ref{string}).

\end{document}